\documentclass[11pt]{article} \usepackage[table]{xcolor}
\usepackage{graphicx} \usepackage{amsmath,amsfonts,amssymb}
\usepackage{graphicx} \usepackage{multirow} \usepackage{tikz,pgf}
\usepackage{url} \usepackage{amsmath} \usepackage{graphicx}
\usepackage{epsfig} \usepackage{epstopdf} \usepackage{color}
\usepackage{enumitem} \usepackage{ amsfonts} \usepackage{amsmath}
\usepackage{amsthm} \parindent 0em \parskip 0.5em
 \def\tto{\;{\lower 1pt
\hbox{$\rightarrow$}}\kern -10pt \hbox{\raise 2pt
\hbox{$\rightarrow$}}\;}  
\def\Bar{\overline}  
 \def\epsilon{\varepsilon} 
\def\h{\hfill\Box} \def\R{\Bbb R} 
 \def\N{\Bbb N} \def\ox{\bar{x}}
 \def\oy{\bar{y}}  
  
 \def\span{\mbox{\rm span}}
 \def\ri{\mbox{\rm ri}}
\def\graph{\mbox{\rm gph}} \def\gph{\mbox{\rm gph}}
\def\aff{\mbox{\rm aff}} \def\epi{\mbox{\rm epi}}
 
\def\dom{\mbox{\rm dom}} \def\aff{\mbox{\rm aff}}  \def\sint{\mbox{\rm int}} \def\rint{\mbox{\rm rint}}
 
     \def\cone{\mbox{\rm cone}} \def\rge{\mbox{\rm rge}}
\def\iri{\mbox{\rm iri}} \def\sqri{\mbox{\rm sqri}}

 \def\h{\hfill\square} 
 \def\emp{\emptyset} 
\def\oR{\Bar{\R}}   
   
   \def\Om{\Omega}
 \def\emp{\emptyset} 
\def\oR{\Bar{\R}}

 \def\qri{\mbox{\rm qri}} \def\sqri{\mbox{\rm
sqri}}  \def\qi{\mbox{\rm qi}}
\def\qri{\mbox{\rm qri}}
\setlist[enumerate,1]{itemsep=0.0ex,parsep=0.5ex,label={\rm(\alph*)},leftmargin=*,
align=left} \newcounter{lk}

\usepackage[colorlinks=true]{hyperref}

\begin{document} \begin{center}
{\sc\bf Generalized Relative Interiors and Generalized Convexity\\ in Infinite Dimensions}\\[1ex]
{\sc Vo Si Trong Long}\footnote{Faculty of Mathematics and Computer
Science, University of Science, Ho Chi Minh City,
Vietnam.}$^,$\footnote{Vietnam National University, Ho Chi Minh
City, Vietnam (email: vstlong@hcmus.edu.vn)}, {\sc Boris S.
Mordukhovich}\footnote{Department of Mathematics,
 Wayne State University,
 Detroit, Michigan 48202, USA (boris@math.wayne.edu). Research of this author was partly supported by the USA National Science Foundation under grants DMS-1808978 and DMS-2204519,
 by the Australian Research Council under grant DP-190100555, and by Project~111 of China under grant D21024.}, {\sc Nguyen Mau  Nam}\footnote{Fariborz Maseeh Department of Mathematics and Statistics,
 Portland State University, Portland, OR 97207, USA (mnn3@pdx.edu). Research of this author was partly supported by the USA National Science Foundation under grant DMS-2136228.}\\[2ex]
{\bf In memory of Diethard Pallaschke,\\
an outstanding mathematician and unique individual} \end{center}
\small{\bf Abstract.} This paper focuses on investigating
generalized relative interior notions for sets in locally convex
topological vector spaces with particular attentions to graphs of
set-valued mappings and epigraphs of extended-real-valued functions.
We introduce, study, and utilize a novel notion of {\em quasi-near
convexity} of sets that is an infinite-dimensional extension of the
widely acknowledged notion of near convexity. Quasi-near convexity
is associated with the quasi-relative interior of sets, which is
investigated in the paper together with other generalized relative
interior notions for sets, not necessarily convex. In this way, we
obtain new results on generalized relative interiors for graphs of
set-valued mappings in convexity and generalized convexity settings.\\[1ex]

\newtheorem{theorem}{Theorem}[section]
\newtheorem{proposition}[theorem]{Proposition}
\newtheorem{remark}[theorem]{Remark}
\newtheorem{lemma}[theorem]{Lemma}
\newtheorem{corollary}[theorem]{Corollary}
\theoremstyle{definition}
\newtheorem{definition}[theorem]{Definition}
\newtheorem{example}[theorem]{Example}

\normalsize

\section{Introduction} The concept of {\em relative interior} for
convex sets is highly important in finite dimensions, as it occupies
a pivotal position in convex analysis and its practical applications
to, e.g., convex optimization. Recognizing its fundamental
importance, there have been significant efforts to explore
appropriate notions of generalized relative interior in
infinite-dimensional spaces. The notions of  {\em quasi-interior},
{\em strong quasi-relative interior}, {\em intrinsic relative
interior}, and {\em quasi-relative interior} for convex sets in
infinite dimensions have been well recognized while playing their
own notable roles in various aspects of convex analysis and
optimization; see
\cite{bao,bauschke2011convex,borwein2003notions,borwein1992partially,boct2012regularity,boct2008regularity,nam2022fenchel,
nam2021quasi,ktz,mor,mordukhovich2022convex,ng2003fenchel,pall,zalinescu2002convex}
with the references and discussions therein.

Among the most significant results of finite-dimensional convex
geometry, we mention Rockafellar's theorems on nonemptiness of the
relative interior of a nonempty convex set in $\R^n$ and the
relative interior representation for the graph of a convex
set-valued mapping; see \cite[Theorem~6.8]{r}. Married to convex
separation, the latter theorem lies at the core of the {\em
geometric approach} for generalized differentiation in convex
analysis developed in \cite{nameasy}. This developed approach
provides an easy and unified way to access many important results of
convex analysis, optimization, and their applications.

The importance of the aforementioned relative interior
representation in finite dimensions calls for its extension to the
infinite-dimensional setting by using appropriate generalized
relative interior notions. In particular, it is natural to consider
representations of the quasi-interior, strong quasi-relative
interior, intrinsic relative interior, and quasi-relative interior
for graphs of convex set-valued mappings in infinite dimensions. To
the best of our knowledge, this question has been addressed only for
the case of quasi-relative interior (see
\cite{nam2022fenchel,nam2021quasi,mordukhovich2022convex}), while it
remains open for the other listed cases of quasi-interior, strong
quasi-relative interior, and intrinsic relative interior of convex
graphs.

Another important unsolved issue in this direction is to go {\em
beyond convexity}. Among various notions of generalized convexity
for sets, the so-called {\em near convexity} (known also as ``almost
convexity") seems to be the most natural to consider first. This
notion actually goes back to Minty \cite{Minty1961} in his study of
maximal monotone operators in finite dimensions. It can be
equivalently formulated as the property that the set in question is
situated between a convex set and its closure, with taking into
account that any (nonempty) finite-dimensional convex set has
nonempty relative interior. In \cite{R1970}, Rockafellar extended
Minty's notion and result to smoothly reflexive Banach spaces in
terms of the very differently formulated notion of ``virtual
convexity" with showing that the latter reduced to \cite{Minty1961}
in finite dimensions. More recently, the near convexity of sets and
associated notion for functions have been consider in
\cite{bmw2013,bkw2008,LM2019,mmw2016,NTY2023} but only in
finite-dimensional spaces.

This paper addresses the aforementioned open questions in convex and
nonconvex settings. To proceed with nonconvex sets in the general
framework of {\em locally convex topological vector} (LCTV) spaces,
we introduce a new notion of {\em quasi-near convexity}, which is an
infinite-dimensional extension of near convexity with the usage of
nonempty quasi-related interior of convex sets in the definition of
the new notion instead of (always nonempty) relative interior of
convex sets in the finite-dimensional near convexity. Note that any
nonempty convex set has nonempty quasi-relative interior in the case
of {\em separable Banach} spaces; see \cite{borwein2003notions}.
Then we establish generalizations of Rockafellar's relative interior
representation theorem for set-valued mappings with quasi-nearly
convex graphs.

Our paper is structured as follows. Section~2 contains some basic
notation and definitions of convex analysis broadly used in the
subsequent material. Section~3 focuses on revisiting a number of
important notions of generalized relative interior with further
clarifications. Section~4 is devoted to the study of the intrinsic
relative interior and strong quasi-relative interior of convex
graphs. Section~5 introduces and investigates a new notion of
quasi-near convexity for nonconvex sets. The final Section~6
provides estimates and representations of the quasi-relative
interior for graphs of quasi-nearly convex set-valued mappings.

\section{Preliminaries} Throughout the paper, we use standard
definitions and notation, which can be founded, e.g., in
\cite{mordukhovich2022convex}. Unless otherwise stated, all the
spaces under consideration are real {\em LCTV} . Recall that the
topological dual of $X$ is denoted by $X^*$ with the canonical
pairing $\left\langle x^*, x\right\rangle:=x^*(x)$ for $x \in X$ and
$x^* \in X^*$. A nonempty set $\Omega\subset X$ is a {\em cone} if
$\lambda w \in \Omega$ for all $w \in \Omega$ and $\lambda \geq 0$.
The {\em closure, conic hull, affine hull}, and {\em linear hull} of
$\Om$ are denoted by $\Bar{\Om}$, $\cone(\Om)$, $\aff(\Om)$, and
$\span (\Om)$, respectively.

Let $f\colon X \rightarrow \overline{\mathbb
R}:=\R\cup\{\pm\infty\}$ be an extended-real-valued function. The
{\em effective domain} and {\em epigraph} of $f$ are defined,
respectively, by \begin{equation*} \dom (f):=\left\{x \in X \mid
f(x)<\infty\right\}\;\mbox{ and }\; \epi(f)=\left\{(x, \lambda) \in
X \times \mathbb{R} \mid f(x) \leq \lambda\right\}.
 \end{equation*}
The function $f$ is called {\em proper} if $\dom f \neq \emptyset$
and $f(x)>-\infty$ for all $x \in X$. We also say that $f$ is {\em
convex} if $\epi(f)$ is a convex set.

Given a set-valued mapping $F\colon X \rightrightarrows Y$, define
the {\em domain}, {\em range}, and {\em graph} of $F$ by
\begin{eqnarray*}
\begin{array}{ll}
&\dom( F):=\left\{x \in X \mid F(x) \neq
\emptyset\right\},\;\rge (F):=\bigcup_{x \in X} F(x),\\
&\gph(F):=\left\{(x, y) \in X \times Y \mid y \in F(x)\right\}.
\end{array}
\end{eqnarray*}
If $\gph(F)$ is a convex set in $X \times Y$, then we say that the {\em mapping} $F$ is {\em convex}.

For a function $f\colon X \rightarrow \overline{\mathbb{R}}$, define
the {\em epigraphical mapping} $E_f\colon X \rightrightarrows
\mathbb{R}$ by \begin{equation} \label{epigraphical}
E_f(x):=\{\alpha \in \mathbb{R} \mid f(x) \leq \alpha\}, \ \; x \in
X. \end{equation} It is easy to verify the equalities $$ \dom
(E_f)=\dom (f)\;\mbox{ and }\;\gph (E_f)=\epi (f). $$ We also
consider the {\em epigraphical range} of $f$ given by $\rge
(f):=\rge (E_f)$.

Given a subset $\Omega$ of $X$, define its {\em polar} by $$
\Omega^{\circ}:=\left\{x^* \in X^* \mid\left\langle x^*,
w\right\rangle \leq 1\;\mbox{ for all }\;w \in\Omega\right\}. $$
Therefore, for a subset $\Theta$ of $X^*$ we have $$
\Theta^{\circ}=\left\{x \in X \mid\left\langle z^*, x\right\rangle
\leq 1\;\mbox{ whenever }\;z^* \in \Theta\right\}. $$ It follows
from the definition that if $\Omega$ is a cone in $X$, then $$
\Omega^{\circ}=\left\{x^* \in X^* \mid\left\langle x^*,
w\right\rangle \leq 0 \text { for all } w \in \Omega\right\}, $$ and
if $\Theta$ is a cone in $X^*$, then $$ \Theta^{\circ}=\left\{x \in
X \mid\left\langle z^*, x\right\rangle \leq 0\;\mbox{ for all }\;z^*
\in \Theta\right\}.
$$

Finally in this section,  we recall the fundamental {\em separation
properties} of sets, which are studied and applied in the subsequent
sections. \begin{definition}\label{defseparate} {\rm Let $\Omega_1$
and $\Omega_2$ be two nonempty subsets of $X$. We say that
$\Omega_1$ and $\Omega_2$ can be {\sc separated} by a closed
hyperplane if there exists $x^* \in X^*\setminus\{0\}$ such that
\begin{equation}\label{key6} \sup \left\{\left\langle x^*,
y\right\rangle \mid y \in \Omega_2\right\} \leq \inf
\left\{\left\langle x^*, x\right\rangle \mid x \in \Omega_1\right\}.
\end{equation} If it holds in addition that
\begin{equation}\label{key7} \inf \left\{\left\langle x^*,
y\right\rangle \mid y \in \Omega_2\right\}<\sup \left\{\left\langle
x^*, x\right\rangle \mid x \in \Omega_1\right\}, \end{equation} then
we say that $\Omega_1$ and $\Omega_2$ can be {\sc properly
separated} by a closed hyperplane.} \end{definition}

Observe that (\ref{key6}) can be rewritten as $$ \left\langle x^*,
y\right\rangle \leq\left\langle x^*, x\right\rangle \text { whenever
} y \in \Omega_2 \text { and } x \in \Omega_1, $$ while (\ref{key7})
means that there exist $\bar{y} \in \Omega_2$ and $\bar{x} \in
\Omega_1$ satisfying $$ \left\langle x^*,
\bar{y}\right\rangle<\left\langle x^*, \bar{x}\right\rangle.
$$

\section{Extended Relative Interiors of Sets}

In this section, we first revisit the major {\em generalized
relative interior} properties of sets in LCTV spaces that are
broadly used in the literature; see, e.g.,
\cite{bao,bauschke2011convex,
borwein2003notions,borwein1992partially,cammaroto2005separation,mor,
mordukhovich2022convex,ng2003fenchel,zalinescu2002convex,zalinescu2015}.
We also present here some refinements of known results in the case
of arbitrary (not necessarily convex) sets.

\begin{definition} \label{def1} {\rm Let $\Omega$ be a  subset of
$X$. \begin{enumerate} \item The {\sc interior} of $\Omega$ {\em
with respect to the affine hull} $\aff (\Om)$ is the set $$
\operatorname{rint}(\Omega):=\{x \in \Omega \mid \exists \text{ a
neighborhood } V \text{ of the origin } \text{such that }(x+ V)\cap
\aff(\Omega) \subset \Omega\}. $$ \item  The {\sc interior} of
$\Omega$ {\em with respect to the closed affine hull}
$\overline{\aff} (\Om)$ is the set $$ \operatorname{ri}(\Omega):=\{x
\in \Omega \mid \exists \text{ a  neighborhood } V \text{ of the
origin } \text{such that } (x+ V)\cap \overline{\aff}(\Omega) \subset
\Omega\}. $$ \item The {\sc quasi-interior} of $\Omega$ is the set
$$ \operatorname{qi}(\Omega):=\{x \in \Omega \mid
\overline{\operatorname{cone}}(\Omega-x) =X\}.$$ \item The {\sc
strong quasi-relative interior} of $\Omega$ is the set $$
\operatorname{sqri}(\Omega):=\{x \in \Omega \mid
\operatorname{cone}(\Omega-x) \text { is a closed subspace of } X\}
. $$ \item  The {\sc intrinsic relative interior} of $\Omega$ is the
set $$\operatorname{iri}(\Omega):=\{x \in \Omega \mid
\operatorname{cone}(\Omega-x) \text { is a subspace of } X\}. $$
\item  The {\sc quasi-relative interior} of $\Omega$ is the set $$
\operatorname{qri}(\Omega):=\{x \in \Omega \mid
\overline{\operatorname{cone}}(\Omega-x) \text { is a subspace of }
X\}. $$ \end{enumerate}}
If $\qri(\Om)=\iri(\Om)$, we say that $\Om$ is {\sc QUASI-REGULAR.}
\end{definition}

The following proposition establishes a relationship between
$\mbox{\ri}(\Omega)$ and $\mbox{\rm rint}(\Omega)$.

\begin{proposition}\label{rirep} Let $\Omega$ be a subset of $X$,
not necessarily convex. Then we have
\begin{equation}\label{rirepf}
\ri(\Omega)=\begin{cases}\mbox{\rm rint}(\Omega)& \mbox{ if }\aff(\Omega)\; \mbox{ is closed},\\
\emptyset & \mbox{ otherwise}. \end{cases} \end{equation}
In particular, $\mbox{\rm ri}(\Omega)\subset\mbox{\rm rint}(\Omega)$,
where the equality holds if ${\rm ri}(\Omega)\neq\emptyset$.
\end{proposition} \noindent {\bf Proof}. If $\aff(\Omega)$ is
closed, then by definition $\ri(\Omega)=\mbox{\rm rint}(\Omega)$.
Consider the case where $\aff(\Omega)$ is not closed and suppose on
the contrary that $\ri(\Omega)\neq\emptyset$. Pick $\ox\in
\ri(\Omega)$ and find a neighborhood $V$ of the origin such that
\begin{equation}\label{ridef} (\ox+V)\cap
\overline{\aff}(\Omega)\subset \Omega. \end{equation} Fix any $z\in
\overline{\aff}(\Omega)$ and select $\lambda>0$ so small that
$\ox+\lambda(z-\ox)\in \ox+V$. Then we get
$\oy:=\ox+\lambda(z-\ox)=\lambda z+(1-\lambda)\ox\in (\ox+V)\cap
\overline{\aff}(\Omega)$ because $\oy$ is an affine combination of
$z,\ox\in \overline{\aff}(\Omega)$. Using \eqref{ridef} gives us
$\oy\in \Omega$. Then \begin{equation*}
z=\left(1-\frac{1}{\lambda}\right)\ox+\frac{1}{\lambda}\oy\in
\aff(\Omega) \end{equation*} since $z$ is an affine combination of
$\ox, \oy\in \Omega$. This yields the closedness of $\aff(\Omega)$,
a contradiction completing the proof of \eqref{rirepf}. The last
statement follows from \eqref{rirepf}. $\h$

The next result provides relationships between the notions of
generalized relative interiors from Definition~\ref{def1}.

\begin{theorem}\label{2s} Let $\Omega$ be a subset of $X$, not
necessarily convex. Then we have the inclusions
\begin{equation}\label{s1} \ri(\Omega)\subset\mbox{\rm
sqri}(\Omega)\subset \iri(\Omega)\subset \mbox{\rm qri}(\Omega),
\end{equation} \begin{equation}\label{s2} \ri(\Omega)\subset
\mbox{\rm rint}(\Omega)\subset \iri(\Omega). \end{equation} All the
inclusions above become equalities if $\Omega$ is convex with
$\ri(\Omega)\neq\emptyset$. \end{theorem} \noindent {\bf Proof}.
Take any $\ox\in \ri(\Omega)$. Then $\ox\in \Omega$ and there exists
a neighborhood $V$ of the origin such that \eqref{ridef} is
satisfied. Let us verify that \begin{equation}\label{affcone}
\overline{\aff}(\Omega)-\ox =\cone(\Omega-\ox). \end{equation}
Indeed, take any $z\in  \overline{\aff}(\Omega)-\ox$ and find
$\lambda>0$ so small that $\lambda z\in V$. Since
$\overline{\aff}(\Omega)-\ox$ is a linear subspace, we have $\lambda
z\in \overline{\aff}(\Omega)-\ox$ and thus \begin{equation*}
\ox+\lambda z\in (\ox+V)\cap \overline{\aff}(\Omega)\subset \Omega.
\end{equation*} Then $z\in 1/\lambda(\Omega-\ox)\subset
\cone(\Omega-\ox)$, which justifies the inclusion $``\subset"$ in
\eqref{affcone}. Observe that $\Omega-\ox\subset
\overline{\aff}(\Omega)-\ox$. Since the latter set is a linear
subspace, we arrive at the reverse inclusion in \eqref{affcone}.

By \eqref{affcone}, $\cone(\Omega-\ox)$ is a closed linear subspace
of $X$ and thus $\ox\in \sqri(\Omega)$, which yields
$\ri(\Omega)\subset \sqri(\Omega)$. The inclusion $\mbox{\rm
sqri}(\Omega)\subset \iri(\Omega)$ follows from the obvious fact
that any closed linear subspace is a linear subspace. Since the
closure of a linear subspace is also a linear subspace, we obtain
the inclusion $\iri(\Omega)\subset \mbox{\rm qri}(\Omega)$ and thus
complete the proof of \eqref{s1}.

The inclusion $\ri(\Omega)\subset \mbox{\rm rint}(\Omega)$ is a
consequence of Proposition~\ref{rirep}. Taking now any $\ox\in
\mbox{\rm rint}(\Omega)$ and following the proof of \eqref{affcone}
tell us that \begin{equation*} \aff(\Omega)-\ox =\cone(\Omega-\ox),
\end{equation*} which implies that $\cone(\Omega-\ox)$ is a linear
subspace. Hence $\ox\in \iri(\Omega)$ verifying \eqref{s2}.

Finally, assume that $\Omega$ is convex with
$\ri(\Omega)\neq\emptyset$. Then $\ri(\Omega)=\qri(\Omega)$ (see
\cite{borwein2003notions} and \cite{mordukhovich2022convex}), which
clearly implies that all the inclusions in \eqref{s1} and \eqref{s2}
become equalities. $\h$

Next, we provide an example to demonstrate that the inclusion $\rint(\Om)\subset \iri(\Om)$ is strict in general.
\begin{example}
Consider the set $\Om:=\{(x,\lambda)\in\mathbb R^2\mid x^2\leq \lambda\}\cup \big(\mathbb R\times (-\infty,0]\big)$. Then we have
$\cone (\Om-(0,0))=\mathbb R^2$, which implies that $(0,0)\in \iri(\Om)$. However, it is easy to verify that $(0,0)\notin \rint(\Om).$
\end{example}

The example below shows that $\text{\rm rint}(\Omega)$ and $\sqri(\Om)$ are different a convex set $\Om$. More examples that distinguish other notions of generalized relative interiors can be found in \cite{bauschke2011convex,borwein1992partially}.

\begin{example}\label{c01}{\rm  Let $X:=C_{[0,1]}$ (the normed space of real continuous functions on $[0,1]$) with the ``max'' norm, and let $P$ be the
set of all polynomials with real coefficients on $[0,1]$. It is well known that $P$ is
a dense subspace in $X$. Thus we get $$ \rint(P)=P\neq
\emptyset=\sqri(P). $$}
\end{example}

To proceed further, take an affine subset $M$ of an LCTV space $X$.
It is well known that there exists a unique linear subspace $L$ of
$X$ such that \begin{equation*} M=x_0+L\; \mbox{\rm for some
}\;x_0\in M. \end{equation*} In this case, $L$ is called the {\em
linear subspace parallel to} $M$. We have the representation
$L=M-M$; see, e.g., \cite{mordukhovich2022convex}. Consider now a
nonempty subset $\Omega$ of $X$ and get that $M:=\aff(\Omega)$ is a
nonempty affine set in $X$. The next result provides a
representation of the linear subspace that is parallel to $\aff
(\Omega)$.

\begin{proposition} Let $\Omega$ be a nonempty set in $X$. Then
$L:=\mbox{\rm aff}(\Omega-\Omega)$ is the linear subspace which is
parallel to $\mbox{\rm aff}(\Omega)$. \end{proposition} {\bf
Proof.} By the above, the linear subspace parallel to $\aff(\Omega)$
is $\aff(\Omega)-\aff(\Omega)$. It is an easy exercise to show that
$\aff (\Omega)-\aff (\Omega)=\mbox{\rm
aff}(\Omega-\Omega)=L$\footnote{In general, for two sets $A,B\subset
X$ and numbers $\alpha,\beta\in \R$ we have $$\mbox{\rm aff}(\alpha
A+\beta B)=\alpha\, \aff (A)+\beta\,\aff (B).$$ {\em Hint}: Pick
$x_0\in A$ and $y_0\in B$. It suffices to show that
\begin{equation*} \mbox{\rm aff}(\alpha C+\beta D)=\alpha \,\aff
(C)+\beta \,\aff( D), \end{equation*} where $C:=A-x_0$ and
$D:=B-x_0$, which both contain the origin.}. This completes the
proof. $\h$

The following lemma presents a straightforward result involving the convexity of the conic hull of a set $\Omega$ given by $\cone(\Omega):=\{\lambda x\; |\; x\in \Omega\}$.

\begin{lemma}\label{l2} If $\Omega$ is a convex set in $X$,
then so is $\cone\, (\Omega)$.
 \end{lemma} {\bf Proof.} Take
any $x, y\in \cone\,(\Omega)$. Then $x=\lambda_1w_1$ and
$y=\lambda_2w_2$, where $\lambda_1, \lambda_2\geq 0$ and $w_1,
w_2\in \Omega$. Assuming that $\lambda_1+\lambda_2>0$, we get
$$x+y=\lambda_1w_1+\lambda_2w_2=(\lambda_1+\lambda_2)\left(\frac{\lambda_1}{\lambda_1+\lambda_2}
w_1+\frac{\lambda_2}{\lambda_1+\lambda_2}w_2\right)\in
(\lambda_1+\lambda_2)\Omega\subset \cone\,(\Omega).$$ This yields
the convexity of the set $\cone(\Omega)$ since it is a cone. $\h$

Although some results of the  next three propositions can be
distilled from \cite{zalinescu2002convex}, we prefer for the
reader's convenience and completeness to give here their simplified
proofs with commentaries.

\begin{proposition}\label{prss} If $\Om$ is a nonempty convex set in $X$,
then \begin{equation*} \mbox{\rm span} (\Om-w)=\cone(\Om-\Om)
\end{equation*} whenever $w\in \Om$\footnote{The conclusion does not
hold in general. Indeed, we can consider $\Omega=\{(1,0)\cup (0,
1)\}\subset\R^2$.}. \end{proposition} {\bf Proof.} Since $\Om-\Om$
is convex, by Lemma~\ref{l2} the set $\cone(\Om-\Om)$ is a
convex cone. It is indeed a linear subspace because if $x\in
\cone(\Om-\Om)$, then $-x\in \cone(\Om-\Om)$. For any $w\in \Om$, we
have $\Om-w\subset \cone(\Om-\Om)$ and thus $\mbox{\rm
span}(\Om-w)\subset\cone(\Om-\Om)$.

To verify the reverse inclusion, fix any $x\in \cone(\Om-\Om)$ and
get the representation \begin{equation*} x=\lambda(w_1-w_2)\;
\mbox{\rm where }\lambda\geq 0,\text{ and } w_1, w_2\in \Om. \end{equation*}
Then $x=\lambda(w_1-w)-\lambda(w_2-w)\in \mbox{\rm span}(\Om-w)$.
This justifies the reverse inclusion and completes the proof. $\h$

\begin{proposition}\label{affcone1} Let $\Om$ be a nonempty convex
set in $X$, and let $w\in \Om$. Then we have: \begin{enumerate} \item
$\aff(\Om)-w=\cone(\Om-\Om)$. \item $\cone(\Om-\Om)$ is the linear
subspace parallel to $\aff(\Om)$. \item
$\aff(\Omega)-w=\cone(\Omega-w)$ if and only if $\cone(\Om-w)$ is a
linear subspace of $X$. \item
$\overline{\aff}(\Om)-w=\overline{\cone}(\Om-w)$ if and only if
$\overline{\cone}(\Om-w)$ is a linear subspace of $X$.
\end{enumerate} \end{proposition} {\bf Proof.} (a) It follows from
Proposition~\ref{prss} that \begin{equation*}
\cone(\Om-\Om)=\mbox{\rm span}(\Om-w). \end{equation*} Since the set
$\aff(\Om-w)$ is a linear subspace that contains $\Om-w$, we see
that \begin{equation*} \cone(\Om-\Om)= \mbox{\rm span}(\Om-w)\subset
\aff(\Om-w)=\aff(\Omega)-w. \end{equation*} Observe that
$\Om-w\subset \Om-\Om\subset \cone(\Om-\Om)$, where the last set is
a linear subspace and hence an affine set. Thus we have the inclusion $\aff(\Om-w)\subset \cone(\Om-\Om)$, which
completes the proof of assertion (a).

(b) This assertion follows directly from (a).

(c) If $\aff(\Om)-w=\cone(\Om-w)$, it is obvious that $\cone(\Om-w)$
is a linear subspace since $\aff(\Om)-w$ has this property. To
verify the converse implication, suppose that $\cone(\Om-w)$ is a
linear subspace. Observe by (a) that \begin{equation*}
\cone(\Om-w)\subset \cone(\Om-\Om)=\aff(\Om)-w. \end{equation*} Take
further any $x\in \aff(\Om)-w=\cone(\Om-\Om)$ and find $\lambda\geq
0$ and $w_1,w_2\in \Om$ such that \begin{equation*}
x=\lambda(w_1-w_2)=\lambda(w_1-w)+\lambda(w-w_2). \end{equation*}
Since $w_2-w\in \cone(\Om-w)$, where the latter set is a linear
subspace, we see that $w-w_2\in \cone(\Om-w)$. This shows that $x\in
\cone(\Om-w)$ and hence justifies (c).

(d) Assume that $\overline{\cone}(\Om-w)$ is a linear
subspace of $X$. It follows from (a) that \begin{equation*}
\overline{\cone}(\Om-w)\subset
\overline{\cone}(\Om-\Om)=\overline{\aff(\Om)-w}=\overline{\aff}(\Om)-w.
\end{equation*} To verify the reverse inclusion, it suffices to show
that $\cone(\Om-\Om)\subset \overline{\cone}(\Om-w)$. Take any $x\in
\cone(\Om-\Om)$ and find $\lambda\geq 0$, $w_1,w_2\in
\Om$ such that \begin{equation*}
x=\lambda(w_1-w_2)=\lambda(w_1-w)+\lambda(w-w_2). \end{equation*}
Similarly to the proof of (c), we see that $w-w_2\in
\overline{\cone}(\Om-w)$, which implies that $x$ belongs to the set
$\overline{\cone}(\Om-w)$ since the latter is a linear subspace.
This completes the proof of (d) by taking into account that the
other implication is obvious. $\h$

The relationship between the strong quasi-relative interior and the
intrinsic relative interior of a convex set is established next.

\begin{proposition}\label{si} Let $\Omega$ be a nonempty convex set in $X$.
Then we have the representation \begin{equation}\label{sqrirep}
\sqri(\Omega)=\begin{cases}\mbox{\rm iri}(\Omega)& \mbox{ if }\;\aff(\Omega)\; \mbox{ is closed},\\
\emptyset & \mbox{ otherwise}. \end{cases} \end{equation}
\end{proposition} \noindent {\bf Proof}. Suppose that $\aff(\Omega)$
is closed. Taking any $\ox\in \iri(\Omega)$, it follows from
definition that $\cone(\Omega-\ox)$ is a linear subspace of $X$. By
Proposition~\ref{affcone1}(c) we have
$\cone(\Om-\ox)=\aff(\Om)-\ox$, and so $\cone(\Om-\ox)$ is a closed
linear subspace of $X$. Thus $\ox\in \sqri(\Omega)$, which implies
that $\iri(\Omega)\subset\sqri(\Omega)$. Since the reverse inclusion
follows from Theorem~\ref{2s}, we justify the equality
$\sqri(\Omega)=\iri(\Omega)$ in this case.

Assume now that $\aff(\Om)$ is not closed. Arguing by contradiction,
suppose that $\sqri(\Om)$ is nonempty and pick $\ox\in \sqri(\Om)$.
Then $\cone(\Om-\ox)$ is a closed linear subspace of $X$. It follows
from Proposition~\ref{affcone1}(c) that $\aff(\Om)-\ox=\cone(\Om-\ox)$ is
also a closed linear subspace of $X$. Hence the affine hull
$\aff(\Om)$ is closed, which is a contradiction verifying that the
set $\sqri(\Om)$ is empty in this case. $\h$

The next proposition shows that the notions of
generalized relative interiors in Definition~\ref{def1} for convex set agree with those taken from \cite{bauschke2011convex} under different
names.

\begin{proposition} \label{pro3.9} Let $\Omega$ be a convex
set in $X$. We have the following assertions: \begin{enumerate}
\item $\sqri(\Om)=\{x\in \Om\; |\; \cone(\Om-x)=\overline{\mbox{\rm
span}}(\Om-x)\}$. \item $\iri(\Om)=\{x\in \Om\; |\;
\cone(\Om-x)=\mbox{\rm span}(\Om-x)\}$. \item $\qri(\Om)=\{x\in
\Om\; |\; \overline{\cone}(\Om-x)=\overline{\mbox{\rm
span}}(\Om-x)\}$. \end{enumerate} \end{proposition} \noindent {\bf
Proof}. (a) Suppose that $x\in \sqri(\Om)$. Then $x\in \Om$ and
$\cone(\Om-x)$ is a closed linear subspace. By
Propositions~\ref{prss} and \ref{affcone1}(a,c) we have
\begin{equation*} \cone(\Om-x)=\aff(\Om)-x=\cone(\Om-\Om)=\mbox{\rm
span}(\Om-x). \end{equation*} It follows that \begin{equation*}
\cone(\Om-x)=\overline{\cone}(\Om-x)=\overline{\mbox{\rm
span}}(\Om-x). \end{equation*}
This justifies the inclusion $``\subset"$ of the equality in (a). The reverse inclusion follows directly from the definition of strong quasi relative interior. \\[1ex]
(b) Suppose that $x\in\iri(\Om)$. Then by definition $\cone(\Om-x)$
is a linear subspace of $X$. Using Propositions~\ref{prss} and
\ref{affcone1}(a,c) as in the proof of assertion (a), we get
\begin{equation*} \cone(\Om-x)=\mbox{\rm span}(\Om-x).
\end{equation*}
This implies the inclusion $``\subset"$ in (b). The reverse inclusion is obvious.\\[1ex]
(c) Assume finally that $x\in \qri(\Om)$ and get that
$\overline{\mbox{\rm cone}}(\Om-x)$ is a linear subspace. Using
again Propositions~\ref{prss} and \ref{affcone1}(a,d) gives us
\begin{equation*}
\overline{\cone}(\Om-x)=\overline{\aff}(\Om)-x=\overline{\cone}(\Om-\Om)=\overline{\mbox{\rm
span}}(\Om-x). \end{equation*} This justifies the inclusion
$``\subset"$ in (c). The reverse inclusion follows from the
definition of quasi-relative interior and the fact that
$\overline{\mbox{\rm span}}(\Om-x)$ is a linear subspace. $\h$

The following result, which extends
\cite[Lemma~2.3]{borwein2003notions} from convex sets to arbitrary
sets in LCTV spaces, plays an important role in representations of
generalized relative interiors for graphs of set-valued mappings.

\begin{proposition} \label{pro1} Let $\Omega$ be a subset of $X$,
and let $\bar x\in\Omega$. Then $\bar x\in\iri(\Omega)$ if and only
if for each $x\in \Omega$ there exists $x_0\in
\Omega$ such that $\bar x=(1- t_0) x+t_0x_0$ with some
$t_0\in(0,1).$ \end{proposition} {\bf Proof}. Suppose $\bar x\in\iri
(\Omega)$. Fix any $x\in \Omega.$ If $x=\bar x$, we immediately get the conclusion. Otherwise, since $\cone(\Omega-\bar x)$ is a linear
subspace, $\bar x-x\in \cone(\Omega-\bar x)$. Thus there exist
$\lambda_0> 0$ and $x_0\in \Omega$ such that $\bar
x-x=\lambda_0(x_0-\bar x)$, which reads as $\bar
x=\frac{1}{1+\lambda_0}x+\frac{\lambda_0}{1+\lambda_0}x_0$. Setting
$t_0:=\frac{\lambda_0}{1+\lambda_0}\in (0,1)$, we obtain $\bar
x=(1-t_0)x+t_0x_0.$

To verify the converse implication, pick any nonzero vector $v\in
\cone(\Omega -\bar x)$. Then we find $\lambda_0>0$ and $x\in
\Omega$ such that $v=\lambda_0(x-\bar x)$. By the
assumption, there exists $x_0\in \Omega$ with $\bar x=(1-t_0)
x+t_0x_0$ for some $t_0\in(0,1)$. It follows that
$$-v=\lambda_0\left(\bar x-\frac{\bar x-t_0 x_0}{1-t_0 }\right)=
\frac{\lambda_0t_0}{1-t_0}\left(x_0-\bar x\right). $$ Consequently,
$-v\in \cone(\Omega -\bar x)$, which completes the proof. $\h$

\begin{proposition}\label{pro3.12} Let $\Om$ be a nonempty convex
set. If $\bar x\in \mbox{\rm sqri}(\Om)$ $(\ri(\Om)$, $\iri(\Om)$,
$\qri(\Om)$, respectively$)$ and $x_0\in \Om$, then for every $t\in
[0,1)$ we have the inclusions $(1-t)\bar x+tx_0 \in \mbox{\rm
sqri}(\Om)$ $(\ri(\Om)$, $\iri(\Om)$, $\qri(\Om)$, respectively$)$.
\end{proposition} {\bf Proof.} We only provide the proof for the
case of $\sqri(\Om)$. The proofs for the other cases can be found in
\cite{borwein2003notions,mordukhovich2022convex}. Let $w_t=(1-t)\bar
x+tx_0$ for $0\leq t<1$. Then \begin{equation*}
\Om-\Om\supset\Om-w_t=\Om-(1-t)\bar x-tx_0=(1-t)(\Om-\bar
x)+t(\Om-x_0)\supset (1-t)(\Om-\bar x). \end{equation*} It follows
therefore that $$\Bar{\cone}(\Om-\Om)\supset \cone (\Om-w_t)\supset
\cone(\Omega-\bar x).$$ As seen in the proof of
Proposition~\ref{pro3.9} (a), we have \begin{equation*}
\overline{\cone}(\Om-\Om)=\cone(\Om-w_t)=\cone(\Omega-\bar x),
\end{equation*} which implies by definition that $w_t\in \mbox{\rm sqri}\,(\Om)$
and thus completes the proof of the proposition. $\h$

\begin{proposition} \label{linearmapiri} Let $T\colon X\to Y$ be a
linear mapping, and let $\Omega$ be a subset of $X$. Then $$T(\iri
(\Omega))\subset \iri (T(\Omega)),$$ where the equality holds if
$\Om$ is convex and $\iri(\Om)\neq\emptyset$. \end{proposition}
\noindent \textbf{Proof.} Fixing any $\bar x \in \iri(\Omega)$, we
get that $\cone(\Omega -\bar x)$ is a linear subspace of $X$. It
follows from the linearity of $T$ that the set
$\cone(T(\Omega)-T(\bar x)) = T(\cone(\Omega - \bar x))$ is a linear
subspace of $Y$. Thus $T(\bar x) \in \iri(T(\Omega))$, which
justifies the claimed inclusion.

Now assume that $\iri(\Omega)\neq\emptyset$ and fix $\bar{x}\in
\iri(\Omega)$. Letting $\bar{y}=T(\bar{x})$, observe that
$\bar{y}\in \iri(T(\Omega))$. Picking any $y\in
\iri(T(\Om)$, we deduce from
Proposition~\ref{pro1} that there exists $y_0\in T(\Om)$ such that
$y= (1-t_0)\bar{y}+ t_0y_0$ for some $t_0\in(0,1)$. Choose $x_0\in
\Om$ with $y_0=T(x_0)$. Then \begin{equation*} y=(1-t_0) \bar{y}+t_0
y_0=(1-t_0) T(\bar{x})+t_0T(x_0)=T((1-t_0)\bar{x}+t_0x_0).
\end{equation*} Letting further $x:=(1-t_0)\bar{x}+t_0x_0$ and
taking into account the assumed convexity of $\Om$ allow us to apply
Proposition \ref{pro3.12} and conclude that $x\in \iri(\Om)$. Thus
$y\in T(\iri(\Om))$, which yields the reverse inclusion and hence
completes the proof of the proposition. $\h$

To study below images of strong quasi-relative interiors of sets
under linear mappings, we need the following definition.

\begin{definition} {\rm A linear mapping $T\colon X\to Y$ is said to
be {\sc subspace closed} if it maps any closed subspace of $X$ to a
closed subspace of $Y$}. \end{definition}

First, we present an easy verifiable sufficient condition for subspace
closedness.

\begin{proposition} Let $T\colon X\to Y$ be a continuous linear mapping,
and let $X$ be complete. Suppose that there exists $\gamma>0$ such
that \begin{equation*} \gamma\|x\|\leq \|T(x)\|\; \mbox{ for all
}\;x\in X. \end{equation*} Then the mapping $T$ is subspace closed.
\end{proposition} {\bf Proof}. Take any convergent sequence $\{y_n\}\subset T(Z)$, where $Z$
is a closed subspace of $X$. Then $y_n=T(x_n)$ for $x_n\in Z$,
$n\in\N$. We have \begin{equation*} \gamma\|x_m-x_n\|\leq
\|T(x_m)-T(x_n)\|=\|y_m-y_n\|\to 0\; \mbox{ as }\;m,n\to \infty.
\end{equation*} This yields $\{x_n\}$ is a Cauchy sequence, and hence it
converges to $x_0\in Z$ due to the completeness of $X$. Then
$y_n=T(x_n)\to T(x_0)\in T(Z)$ as $n\to \infty$, which verifies that
$T(Z)$ is closed. $\h$

In the next proposition, we use the notion of {\sc strong
quasi-regular} for a nonempty set $\Om\subset X$ meaning that
$\sqri(\Om)=\iri(\Om)$.

\begin{proposition} Let $T\colon X\to Y$ be a continuous linear mapping which is subspace closed, and let $\Om$ be a nonempty convex set in $X$.
Then we have \begin{equation*} T(\mbox{\rm sqri}(\Om))\subset
\mbox{\rm sqri}(T(\Om)). \end{equation*} If in addition
$\sqri(\Om)\neq\emptyset$ and $\Om$ is strongly quasi-regular, then
the reverse inclusion holds. \end{proposition} {\bf Proof}. Suppose
that $w\in \mbox{\rm sqri}(\Om)$. Then by Proposition~\ref{pro3.9}(a) we have
\begin{equation*}
\cone(\Om-w)=\overline{\mbox{\rm span}}\,(\Om-w). \end{equation*} Therefore, we obtain the equalities \begin{equation*}
\cone(T(\Om)-T(w))=T(\cone(\Om-w))=T(\overline{\mbox{\rm
span}}\,(\Om-w))=\overline{\mbox{\rm span}}(T(\Om)-T(w)),
\end{equation*} which tell us that $T(w)\in \mbox{\rm sqri}(T(\Om))$.
To verify the reverse inclusion, assume that
$\sqri(\Om)\neq\emptyset$ and $\Om$ is strongly quasi-regular. Then
Proposition~\ref{linearmapiri} ensures that \begin{equation*}
\sqri(T(\Om))\subset \iri(T(\Om))=T(\iri(\Om))=T(\sqri(\Om)),
\end{equation*} which thus completes the proof of the claimed
result. $\h$

\section{Generalized Relative Interiors of Graphical Sets}

This section is devoted to extending Rockafellar's relative interior
representation theorem in $\R^n$ to generalized relative interiors
of graphs of set-valued mappings and epigraphs of
extended-real-valued functions.

We start with considering {\em intrinsic relative interiors} for
graphs of set-valued mappings.

\begin{theorem} \label{thm_iri_graph} Let $F\colon X
\rightrightarrows Y$ be a set-valued mapping. Then
\begin{equation}\label{irigpheq}
\iri(\gph(F))
\subset \left\{(x, y) \in X \times Y \mid x \in \iri(\dom(F)), y \in
\iri(F(x))\right\}.
\end{equation}
Assuming in addition that $F$ is convex, we have the equality $$\iri(\gph(F)) = \left\{(x, y) \in X
\times Y \mid x \in \iri(\dom(F)), y \in \iri(F(x))\right\}. $$
\end{theorem}
\textbf{Proof.} To verify the inclusion ``$\subset$",
take any $(\bar{x}, \bar{y}) \in \iri(\gph(F))$. Considering the {\em projection mapping $\mathcal{P}\colon X\times Y\to X$ given by
\begin{equation*}
    \mathcal{P}(x, y)=x, \; (x, y)\in X\times Y,
\end{equation*}
gives us $\mathcal{P}(\gph(F))=\dom(F)$. Using Proposition~\ref{linearmapiri}, we have
\begin{equation}\label{equaliri}
    \mathcal{P}(\iri(\gph(F)))\subset \iri(\dom(F)),
\end{equation}
which implies that $\ox\in \iri(\dom(F))$}. To show next that $\bar{y} \in \iri(F(\bar{x}))$,
fix any $y \in F(\bar{x})$ telling us that
$(\bar{x}, y) \in \gph(F)$. By
Proposition~\ref{pro1}, find $(x_0, y_0) \in \gph(F)$ such that $$
(\bar{x}, \bar{y})=(1-t_0)(\bar{x}, y)+t_0(x_0, y_0)\;\text{ for
some }\;t_0\in (0,1). $$ This yields $x_0=\bar{x}$ and $\bar{y}=(1-t_0)
y+t_0 y_0$ with $y_0 \in F(\bar{x})$. It follows therefore by Proposition~\ref{pro1} that
$\bar{y} \in \iri(F(\bar{x}))$, which completes the proof of the
inclusion ``$\subset$".

To verify now the reverse inclusion under the convexity of $\gph
(F)$, fix $\bar{x} \in \iri(\dom (F))$ and $\bar{y} \in
\iri(F(\bar{x}))$. Arguing by contradiction, suppose that $(\bar{x},
\bar{y}) \notin\iri (\gph(F))$. Then Proposition~\ref{pro1} yields
the existence of $(x',y') \in \gph (F) $
such that \begin{equation}\label{key3} (\bar x,\bar y)\neq
(1-t)(x',y')+t(x,y)\text{ for all } (x,y)\in \gph (F) \text{ and } t\in (0,1). \end{equation} Consider the following two cases:

\noindent (A) $\bar x\neq (1-t)x'+tx \text{ for all } (x,t)\in \dom
(F)\times (0,1)$,

\noindent (B) $\bar x= (1-t_0)x'+t_0x_0 \text{ for some }
(x_0,t_0)\in \dom (F)\times (0,1)$.

\noindent Let us show below that in each of these cases we arrive at
a contradiction. In case (A),  we clearly have a
contradiction due to $\bar x\in\iri(\dom (F))$.

\noindent In case (B), we distinguish  two subcases:

\noindent (B$_1$) $x' = \bar x$. Since $\bar x= (1-t_0)x'+t_0x_0$,
we get $\bar x= x'=x_0$ and hence $F(\bar x)=F(x')=F(x_0)$. Then it
follows from (\ref{key3}) that $\bar y\neq (1-t)y'+ty $ for all
$(y,t)\in F(\bar x)\times (0,1)$. By Proposition~\ref{pro1} we have
$\bar y\notin\iri(F(\bar x))$, which gives us a contradiction.

\noindent (B$_2$) $x'\neq \bar x$. Note first that since $x_0\in\dom
(F)$, there exists $y_0\in Y$ such that $(x_0,y_0)\in\gph (F)$.
Using the convexity of $\gph (F)$ ensures that
$$(1-t_0)(x',y')+t_0(x_0,y_0)\in \gph(F).$$ Define further the
vector $$\bar y':=(1-t_0)y'+t_0y_0.$$ Since $\bar x=
(1-t_0)x'+t_0x_0$, we get $\bar y'\in F(\bar x)$ and
\begin{equation}\label{key4} (\bar x, \bar
y')=(1-t_0)(x',y')+t_0(x_0,y_0). \end{equation} On the other hand, it
follows from $\bar y\in\iri (F(\bar x))$ and Proposition~\ref{pro1}
that there exists $ y'_0\in F(\bar x)$ such that $\bar y =(1-t'_0)\bar
y'+t'_0y'_0$ for some $t'_0\in (0,1)$. This together with
(\ref{key4}) ensures the equalities $$ \begin{array}{ll}
(\bar x,\bar y)&=(1-t'_0)(\bar x, \bar y')+t'_0(\bar x,y'_0)\\
&=(1-t'_0)[(1-t_0)(x',y')+t_0(x_0,y_0)]+t'_0(\bar x,y'_0)\\
&=(1-t'_0)(1-t_0)(x',y')+(1-t'_0)t_0(x_0,y_0)+t'_0(\bar x,y'_0)\\
&=(1-t'_0)(1-t_0)(x',y')+(t_0-t'_0t_0+t'_0)\left [\frac{(t_0-t'_0t_0)}{(t_0-t'_0t_0+t'_0)}(x_0,y_0)+\frac{t'_0}{(t_0-t'_0t_0+t'_0)}(\bar x,y'_0)\right ]\\
&=(1-s)(x',y')+s(x'',y''), \end{array}$$ where $s:=t_0-t'_0t_0+t'_0\in (0,1)$
and $$
(x'',y''):=\frac{(t_0-t'_0t_0)}{(t_0-t'_0t_0+t'_0)}(x_0,y_0)+\frac{t'_0}{(t_0-t'_0t_0+t'_0)}(\bar x,y'_0)\in \gph (F). $$ This also gives us a contradiction due to \eqref{key3}, which therefore completes the
proof of the theorem. $\h$

\begin{remark}{\rm If in the setting of Theorem~\ref{thm_iri_graph}
an additional assumption $\iri(\gph(F))\ne\emp$ is imposed, then we
have the following {\em alternative simple proof} for the reverse inclusion of \eqref{irigpheq}. Indeed, fix $\bar{x} \in \iri(\dom (F))$ and $\bar{y} \in \iri(F(\bar{x}))$.  By
(\ref{equaliri}) which holds as an equality under the convexity of $\gph(F)$, we find $y_0\in F(\bar x)$ with $(\bar x,y_0)\in
\iri(\gph(F))$. If $y_0=\bar y$, then $(\bar x,\bar y)\in
\iri(\gph(F))$ and we are done. Otherwise, since $\bar y\in
\iri(F(\bar x))$, we deduce from
\cite[Lemma~3.1]{borwein2003notions} that there exists $\gamma >0$
such that $(1+\gamma)\bar y-\gamma y_0\in F(\bar x)$. Then it
follows from the corresponding result of Proposition~\ref{pro3.12}
that $$(1-t)(\bar x,y_0)+t(\bar x, (1+\gamma)\bar y-\gamma y_0)\in
\iri(\gph(F))\;\textrm{ for all }\;t\in [0,1).$$ Choosing
$t:=\frac{1}{1+\gamma}\in (0,1)$, we obtain $$(\bar x, \bar
y)=\left(1-\frac{1}{1+\gamma}\right)\left(\bar
x,y_0\right)+\frac{1}{1+\gamma}\left(\bar x,
(1+\gamma)\bar y-\gamma y_0\right) \in \iri(\gph(F)),$$ which
completes the alternative proof of Theorem~\ref{thm_iri_graph}.}
\end{remark}

To proceed further, we need one more useful result.

\begin{lemma}\label{affgph} Let $F\colon X\tto Y$ be a
set-valued mapping. Suppose that $\mbox{\rm int}(F(x))\neq\emptyset$
for all $x\in \dom(F)$. Then we have the equality \begin{equation*}
\aff(\gph (F))=\aff(\dom (F))\times Y. \end{equation*}
\end{lemma}
 \noindent {\bf Proof}. The inclusion ``$\subset$" is obvious. Let us verify the reverse one. Take any $(x_0, y_0)\in \aff(\dom (F))\times Y$, which yields $x_0\in \aff(\dom (F))$ and $y_0\in Y$. Choose
 $\lambda_i\in \R$ and $x_i\in \dom(F)$ for $i=1, \ldots, m$ such that $\sum_{i=1}^m \lambda_i=1$ and $x_0=\sum_{i=1}^m \lambda_i x_i$. Then select $y_i\in\mbox{\rm int}( F(x_i))$ and define
 $\hat{y}:=\sum_{i=1}^m \lambda_i y_i$ while having in this way that $(x_0, \hat{y})\in \aff(\gph (F))$. We also choose $\epsilon>0$ so small that
 \begin{equation*}
y_i+\epsilon (y_0-\hat{y})\in F(x_i)\; \mbox{\rm for all }i=1,
\ldots, m, \end{equation*} which yields $(x_0, \hat{y}+\epsilon
(y_0-\hat{y}))\in \aff(\gph (F))$, or equivalently
$$(0,\epsilon(y_0-\hat{y}))\in \aff(\gph (F))-(x_0,\hat{y}).$$ Note
that since $(x_0,\hat y)\in \aff(\gph (F))$, the affine hull
$\aff(\gph (F))-(x_0,\hat{y})$ is a linear subspace. This gives us
the inclusion
$$(0,y_0-\hat{y})=\frac{1}{\epsilon}(0,\epsilon(y_0-\hat{y}))\in
\frac{1}{\epsilon}\left ( \aff(\gph (F))-(x_0,\hat{y})\right
)=\aff(\gph (F))-(x_0,\hat{y}),$$ and tells us therefore that
$$(x_0,y_0)=(x_0,\hat y)+(0,y_0-\hat y)\in \aff(\gph
(F))-(x_0,\hat{y})+ (x_0,\hat y)=\aff(\gph (F))$$ and thus completes
the proof of the lemma. $\h$

The next theorem provides a new representation of {\em strong
quasi-relative interiors} for graphs of convex set-valued mappings.

\begin{theorem}\label{sqrigph} Let $F\colon X\tto Y$ be a convex
set-valued mapping. Then we have \begin{equation}\label{sqrif}
\sqri(\gph (F))= \big\{(x, y)\; \big |\; x\in \sqri(\dom (F)), \; y
\in \sqri(F(x)\} \end{equation} provided that $\mbox{\rm
int}(F(x))\neq\emptyset$ for all $x\in\dom (F)$. \end{theorem}
\noindent {\bf Proof}. Take any $(x, y)\in \sqri(\gph (F))$. Then it
follows from Propositions~\ref{si} and \ref{affgph} that the affine
hull $\aff(\gph (F))=\aff(\dom (F))\times Y$ is closed and that
$(x,y)\in \iri(\gph (F))$. By Theorem~\ref{thm_iri_graph}, we have
$x\in \iri(\dom (F))$ and $y\in \iri(F(x))$. Since $\aff(\dom (F))$
is closed, using Proposition~\ref{si} again gives us $x\in
\sqri(\dom F)$. In addition, $y\in \iri(F(x))=\sqri(F(x))$ since
$\mbox{\rm int}(F(x))\neq\emptyset$. This justifies the inclusion
$``\subset"$ in \eqref{sqrif}.

Conversely, take $x\in\sqri(\dom (F))$ and $y\in \sqri(F(x))$. Then
$\aff(\dom (F))$ is closed and $x\in \iri(\dom (F))$. Since
$\mbox{\rm int}(F(x))\neq\emptyset$, we see that $y\in \iri(F(x))$.
It follows from Theorem~\ref{thm_iri_graph} and
Lemma~\ref{affgph} that $(x, y)\in \iri(\gph F)$ and that
$\aff(\gph (F))$ is closed. This yields $(x, y)\in \sqri(\gph (F))$,
which completes the proof. $\h$

The following important statement is a consequence of both
Theorems~\ref{thm_iri_graph} and \ref{sqrigph} applying to {\em
epigraphs} of extended-real-valued convex functions.

\begin{proposition} Let $f\colon X\to \oR$ be an
extended-real-valued convex function. Then we have the generalized
relative interior representations \begin{equation*} \iri(\epi
f)=\big\{(x, \alpha)\in X\times \mathbb R\; \big |\; x\in \iri(\dom
f), \; f(x)<\alpha \big\}, \end{equation*} \begin{equation*}
\sqri(\epi (f))=\big\{(x, \alpha)\in X\times \mathbb R\; \big |\;
x\in \sqri(\dom (f)), \; f(x)<\alpha\big\}. \end{equation*}
\end{proposition} \noindent {\bf Proof}. Consider the epigraphical
mapping $E_f$ associated with $f$ given in \eqref{epigraphical}.
Then the first formula follows from Theorem~\ref{thm_iri_graph}.
Since $\mbox{\rm int}(E_f(x))$ is nonempty for all $x\in \dom
(f)=\dom( E_f)$, the second formula follows from
Theorem~\ref{sqrigph}. $\h$

To conclude this section, we provide a result on the representation
of {\em strong quasi-relative interiors} for {\em ideally convex}
graphs of set-valued mappings, the notion taken from
\cite{holmes2012geometric}.

\begin{definition}\label{defideallyconvex} {\rm A subset $\Omega$ of
$X$ is called {\sc ideally convex} if for any bounded sequence
$\{x_n\} \subset \Omega$ and sequence $\{\lambda_n\}$ of nonnegative
numbers with $\sum_{n=1}^{\infty}\lambda_n = 1$, the series $\sum_{n=1}^{\infty}\lambda_nx_n$ either converges to an element of
$\Omega$, or does not converge at all.} \end{definition}

As observed in \cite{holmes2012geometric}, any convex subset of $X$,
which is either open or closed, is ideally convex. Moreover, every
finite-dimensional convex set is ideally convex.

\begin{proposition} \label{pro_ideallyconvex} Let $F\colon X
\rightrightarrows Y$ be a convex set-valued mapping. Assume that
$\gph (F)$ is an ideally convex subset of the product space $X\times
Y$ and that $\sqri (\gph (F))\neq \emptyset$. Then $$ \sqri(\gph(F))
= \{(x, y) \in X \times Y \mid x \in \sqri(\dom(F)), y \in
\sqri(F(x))\}. $$ \end{proposition} \textbf{Proof.} Combining
Theorem~\ref{thm_iri_graph} with
\cite[Proposition~3.2]{jeyakumar1992generalizations} verifies the
result. $\h$

\section{Generalized Relative Interiors of Quasi-Nearly Convex Sets}

In this section, we introduce a new notion of generalized convexity
called {\em quasi-near convexity}, investigate its basic properties,
and establish its relationship with the near convexity or almost
convexity notions known in the literature; see, e.g., \cite{bmw2013,
bkw2008,LM2019, NTY2023}).

\begin{definition}\label{defqconvex}{\rm Let $\Omega$ be a nonempty subset of $X$. We say that $\Omega$ is  {\sc quasi-nearly convex} if there exists a convex subset $C\subset X$ such that $\qri(C)\neq \emptyset$ and
$$C\subset \Omega \subset \Bar{C}.$$
}\end{definition}

The next proposition shows that the notion of quasi-near convexity agrees with the near convexity in finite dimensions; see, e.g., \cite{bmw2013,bkw2008,LM2019,Minty1961,r,R1970}.

\begin{proposition} \label{equidefq-convex}
Let $\Omega$ be a nonempty set in $\R^n$. The following properties are equivalent:
\begin{enumerate}
\item $\Omega$ is quasi-nearly convex.
\item There exists a convex set $C\subset\mathbb R^n$ such that
\begin{equation*}
\ri(C)\subset\Omega\subset \Bar{C}.
\end{equation*}
\item There exists a convex set $D \subset \mathbb R^n$ such that
\begin{equation*}
D\subset \Omega\subset \Bar{D}.
\end{equation*}
\item  $\Bar{\Omega}$ is convex and $\ri( \Bar{\Omega})\subset\Omega$.
\end{enumerate}
\end{proposition}
\noindent {\bf Proof}. (a)$\Longrightarrow$(b): Obvious. \\[1ex]
(b)$\Longrightarrow$(c): Suppose that (b) is satisfied. Define the convex set $D:=\ri(C)$. Then $\Bar{D}=\Bar{\ri}(C)=\Bar{C}$, and thus
\begin{equation*}
D\subset \Om\subset \Bar{D}.
\end{equation*}
(c)$\Longrightarrow$(d): Suppose that (c) is satisfied. It  can be easily checked that $\Bar{\Om}=\Bar{D}$, and hence $\Bar{\Om}$ is a convex set. In addition, we have
\begin{equation*}
\ri(\Bar{\Om})=\ri(\Bar{D})=\ri(D)\subset D\subset\Om,
\end{equation*}
which ensures that (d) is satisfied.\\[1ex]
(d)$\Longrightarrow$(a): Let (d) hold. Defining $C:=\ri(\Bar{\Om})$, which is nonempty convex set, we get
\begin{equation*}
C=\ri(\Bar{\Om})\subset \Om\subset\Bar{\Om}=\Bar{\ri}(\Bar{\Om})= \Bar{C}.
\end{equation*}
Since $\qri(C)=\ri(C)=\ri(\ri(\Bar{\Om})=\ri(\Bar{\Om})\neq\emptyset$, it follows that (a) is satisfied. $\h$

Some basic properties of quasi-nearly convex sets are established below.

\begin{proposition}\label{pro2} Let $\Omega$ be a quasi-nearly convex set with $C\subset \Omega \subset \Bar{C}$, where $C$ is a convex subset of $X$ satisfying the condition $\qri(C)\neq \emptyset$. Then we have the following:
\begin{enumerate}
\item $\Bar{\Om}=\Bar{C}$.
\item $\qri(C)\subset \qri(\Omega)\subset \qri(\Bar C)$.
\end{enumerate}
Consequently, the set $\qri(\Om)$ is quasi-nearly convex. Moreover, the set $\qri(\Om)$ is convex if $\qri(\Bar{\Om})\subset \Om$, which holds when $X=\R^n$.

\end{proposition}
\textbf{Proof.} (a) The conclusion is an immediate consequence of the definition.

(b) For each $x\in X$, we have
$$
\cone (C-x)\subset \cone (\Omega-x)\subset \cone (\Bar C-x)\subset \overline{\cone} (C-x),
$$
which implies in turn the equalities
$$\overline{\cone} (C-x)=\overline{\cone} (\Om -x)= \overline{\cone} (\overline{C}-x).$$
The latter brings us to the inclusions $$\qri(C)\subset\qri (\Omega)\subset \qri (\overline{C}),$$
which completes the proof of (b).

We therefore have
$$\qri(C)\subset\qri(\Omega)\subset \qri(\overline{C})\subset\overline{\qri}(\overline{C})=\overline{C}=\overline{\qri}(C),$$
where the equalities hold by \cite[Proposition~2.12]{borwein1992partially}. Moreover, we get that $\qri(C)$ is convex by \cite[Lemma~2.9]{borwein1992partially} and that $\qri(\qri(C))=\qri(C)\neq \emptyset$ by \cite[Proposition~2.5(vii)]{boct2008regularity}. Consequently, the set $\qri(\Om)$ is quasi-nearly convex.  If $\qri(\Bar{\Om})\subset \Om$, we can easily check that $\qri(\Bar{C})=\qri(\Bar{\Om})\subset \qri(\Om)$ by definition. This implies that $\qri(\Om)=\qri(\Bar{\Om})=\qri(\Bar{C})$, which is a convex set. This completes the proof of the proposition.$\h$

The following example shows that each inclusion in Proposition~\ref{pro2}(b) can be strict.

\begin{example}\label{ncex1} Let $X$ and $P$ be the given as in Example~\ref{c01}. Consider an element $\ox\in X\setminus P$ and let $\Omega:=P\cup\{\bar x\}$. Then $\Omega$ is {\em not convex} but {\em quasi-nearly convex}, and we have
$$\qri(P)=P\subset \qri(\Omega)=\Omega\subset\qri(\overline{P})=X.$$
\end{example}

Some other useful properties of quasi-nearly convex sets are collected in the next proposition.

\begin{proposition} \label{lem13}
Let $\Omega$ be a quasi-nearly convex set in $X$,  and let $\bar x \in X$. Then we have the following properties:
\begin{enumerate}
\item $\qi (\Omega)=\Omega \cap \qi (\overline{\Omega})\;\text{ and }\;\qri (\Omega)=\Omega \cap \qri (\overline{\Omega})$.
\item $\overline{\qri}(\Omega)=\overline{\Omega}$.
\item $\qri(\Omega)=\qri(\qri(\Omega))$.
\item $\qri (\Omega+\bar x)=\qri (\Omega)+\bar x$.
\end{enumerate}
\end{proposition}
\textbf{Proof.}
(a) It is easy to see that $\qi (\Omega)\subset \Omega \cap \qi (\overline{\Omega})$. Conversely, take any $x\in  \Omega \cap \qi (\overline{\Omega})$. Then  $\overline{\cone}(\Omega -x)=\overline{\cone}(\overline{\Omega}-x)=X$, which means that $x\in \qi (\Omega)$. Similarly, we have $\qri (\Omega)=\Omega \cap \qri (\overline{\Omega})$.

(b) Since $\Om$ is quasi-nearly convex, there exists a convex set $C\subset X$ satisfying $\qri(C)\neq \emptyset$ and $C\subset \Om\subset \Bar{C}$. Then it follows from Proposition~\ref{pro2} and \cite[Proposition~2.12]{borwein1992partially} that $$\Bar{C}=\overline{\qri} (C)\subset \overline{\qri} (\Omega)\subset \overline{\Omega}=\overline{C}.$$  This justifies assertion (b).

(c) It follows directly from (a) and (b) that  $$\qri\left(\qri(\Omega)\right)= \qri(\Omega)\cap\qri(\overline{\qri}(\Omega))=\qri(\Omega)\cap \qri(\overline{\Omega})=\qri(\Omega).$$
(d) Take any $x\in \qri(\Omega) +\bar x$. Then we have $x-\bar x\in \qri(\Omega)$ and hence $\overline{\cone}(\Omega +\bar x-x)$ is subspace, which ensures that $x\in \qri(\Omega +\bar x)$.

Conversely, take any $x\in \qri(\Omega +\bar x)$. Then the set $\overline{\cone}(\Omega +\bar x-x)$ is a subspace, and hence $-(\bar x-x)\in \qri(\Omega)$. This tells us that $x\in \qri (\Omega) +\bar x$ and thus completes the proof. $\h$

To proceed further, we define the {\em normal cone} to a quasi-nearly convex set in the same way as in the case of pure convexity.

\begin{definition}\label{defnormalcone}
{\rm Let $\Omega\subset X$ be a quasi-nearly convex set, and let $\bar x\in\Omega$. The {\sc normal cone} to $\Omega$ at $\bar x$ is
$$N(\bar x; \Omega ): = \{x^*\in X^*\mid \langle x^*,x-\bar x\rangle\leq 0\;\mbox{ for all }\;x\in\Omega\}.$$}
\end{definition}

The next result provides a {\em characterization of quasi-relative interiors} of quasi-nearly convex sets in LCTV spaces.

\begin{proposition}\label{qrisubspace} Let $\Omega\subset X$ be a quasi-nearly convex set, and let $\bar x\in \Omega$. Then $\bar x\in\qri (\Omega)$ if and only if the normal cone $N(\bar x;\Omega)$ is a subspace.
\end{proposition}
\textbf{Proof.}  By the continuity of $x^*\in X^*$ on $X$, we have that $\langle x^*,x-\bar x\rangle\leq 0$ for all $x\in \Omega$ if and only if $\langle x^*,u\rangle\leq 0$ whenever $u\in \overline{\cone}(\Omega -\bar x)$. It follows that $N(\bar x; \Omega)=\overline{\cone}(\Omega -\bar x)^\circ$. This tells us that if $\ox\in \qri(\Om)$, then $\overline{\cone}(\Omega -\bar x)$ is a subspace, and so is $N(\bar x; \Omega)$.

Conversely, we have that the set $\overline{\cone}(\Omega -\bar x)$ is closed and $0\in \overline{\cone}(\Omega -\bar x)$. Furthermore, it follows from Proposition~\ref{pro2} that $\Bar{\Om}$ is convex, and so is $\overline{\cone}(\Bar \Omega -\bar x)=\overline{\cone}(\Omega -\bar x)$. Employing the classical bipolar theorem (see, e.g., \cite[Theorem 1.1.9]{zalinescu2002convex}) gives us the representation
$$N(\bar x; \Omega)^\circ=(\overline{\cone}(\Omega -\bar x)^\circ)^\circ=\overline{\cone}(\Omega -\bar x).$$
This implies that if $N(\bar x; \Omega)$ is a subspace, then so is $\overline{\cone}(\Omega -\bar x)$, and thus $\ox\in \qri(\Om)$. $\h$

The next result, which is often employed in the subsequent sections, presents an equivalent
description of quasi-relative interiors via {\em proper separation} of {\em quasi-nearly convex sets} in LCTV spaces.

\begin{proposition}\label{pro13}
Let $\Omega$ be a quasi-nearly convex set in $X$, and let $\bar x \in \Omega$. Then $\bar x \notin \qri(\Omega)$ if and only if the sets $\{\bar x\}$ and $\Omega$ can be properly separated by a closed hyperplane.
\end{proposition}
\textbf{Proof.}
As seen in Proposition \ref{qrisubspace}, $\bar{x} \in \operatorname{qri}(\Omega)$ if and only if the normal cone $N(\bar{x} ; \Omega)$ is a linear subspace of $X^*$. It follows that $\bar{x}\notin \operatorname{qri}(\Omega)$ if and only if there is $x^*\in N(\bar{x} ; \Omega)$ such that $-x^*\notin N(\bar{x} ; \Omega)$. By the definition of the normal cone, we have
$$\left\langle x^*, x\right\rangle \leq\left\langle x^*, \bar{x}\right\rangle\;\text{ for all }\;x \in \Omega,$$
while $-x^*\notin N(\bar{x} ; \Omega)$ means that
there exists $\hat x\in \Omega$ satisfying
$$ \left\langle x^*, \hat x\right\rangle <\left\langle x^*,\bar x\right\rangle,$$
which completes the proof. $\h$

The example below shows that even if $\Omega$ is {\em convex}, the assumption $\bar x\in \Omega$ in Proposition \ref{pro13} is {\em essential}.
\begin{example}
Let $X$, $P$, and $\ox$ are given as in Example~\ref{ncex1}. We have seen that $\qri(P)=P$ and $\bar x\notin \qri(P)$. Let us now show that the sets $\{\bar x \}$ and $P$ can not be properly separated by a closed hyperplane. Suppose on the contrary that there exists $x^*\in X\setminus\{0\}$ such that
$$\langle x^*,x\rangle\leq \langle x^*,\bar x\rangle\;\text{ for all }\;x\in P,\;\text{ and }\;\langle x^*,\bar y\rangle< \langle x^*,\bar x\rangle\;\text{ for some }\;\bar y\in P.$$
Since $P$ is dense in $X$ and $x^*$ is continuous, we get that $$\langle x^*,z\rangle\leq \langle x^*,\bar x\rangle\;\text{ for all }\;z\in X,$$
which is impossible. Therefore, the sets $\{\bar x\}$ and $P$ cannot be properly separated by a closed hyperplane because the assumption $\bar x\in P$ is violated.
\end{example}

The next theorem establishes {\em proper separation} of two quasi-nearly convex sets.

\begin{theorem}\label{tachqri1}
Let $\Omega_1$ and $\Omega_2$ be quasi-nearly convex subsets of an  LCTV space X. Assume that $\Om_1\cap \Om_2\neq \emptyset$ and that
\begin{equation} \label{quasiregularity}
\qri(\Om_1-\Om_2)=\qri(\Om_1)-\qri(\Om_2).
\end{equation}
Then the sets $\Om_1$ and $\Om_2$ can be properly separated by a closed hyperplane if and only if we have
\begin{equation}\label{qriempty1}
\qri(\Om_1) \cap \qri(\Om_2)=\emptyset.
\end{equation}
\end{theorem}
\textbf{Proof.} The imposed assumptions ensure that $0\in \Om$ and $\qri(\Om)=\qri(\Om_1)- \qri(\Om_2)$, where $\Om:=\Om_1 - \Om_2$. Therefore, if relation (\ref{qriempty1}) holds, then
$$
0 \notin \qri(\Om)=\qri(\Om_1)-\qri(\Om_2).
$$
According to Proposition \ref{pro13}, the sets $\Om_1-\Om_2$ and $\{0\}$ can be properly separated by a closed hyperplane, which clearly ensures the proper separation of the sets $\Om_1$ and $\Om_2$.

Conversely, suppose that $\Om_1$ and $\Om_2$ can be properly separated by a closed hyperplane. Then the sets $\Om_1-\Om_2$ and $\{0\}$ can be properly separated by a closed hyperplane as well. Using Proposition~\ref{pro13} and (\ref{quasiregularity}) yields
$$
0\notin \qri(\Om)=\qri(\Om_1)-\qri(\Om_2),
$$
and hence $\qri(\Om_1) \cap \qri(\Om_2)=\emptyset$, which completes the proof.  $\h$

Based on the sufficient conditions for the fulfillment of (\ref{quasiregularity}) from \cite[Theorem~2.183(c)]{mordukhovich2022convex}, we obtain the following result for the convex case as a direct consequence of Theorem~\ref{tachqri1}; see \cite[Theorem~4.2]{nam2022fenchel} and \cite[Theorem~2.184]{mordukhovich2022convex}.

\begin{corollary}\label{Theoqricap} Let $\Omega_1$ and $\Omega_2$ be convex subsets of  $X$ such that $\Omega_1\cap \Omega_2\neq\emptyset$. Assume that $\qri(\Omega_1)\ne\emp$, $\qri(\Omega_2)\ne\emp$, and the set difference
$\Omega_1-\Omega_2$ is quasi-regular. Then $\Omega_1$ and $\Omega_2$ can be properly separated if and only if \begin{equation*}\qri(\Omega_1)\cap\qri(\Omega_2)=\emp.
\end{equation*}\end{corollary}

To continue our study of quasi-nearly convex sets and their generalized relative interiors, we need the following proposition.

\begin{proposition}\label{qriconvex}
Let $\Omega$ be a quasi-nearly convex set in $X$. If $\bar x \in \qri(\Omega)$, $x_0\in \Omega$, and  $(1-t_0)\bar x+ t_0x_0\in\Omega$ for some $t_0\in (0,1]$, then $(1-t_0)\bar x+ t_0x_0\in\qri (\Omega)$.
\end{proposition}
\textbf{Proof.}
Suppose on the contrary that $\hat x:=(1-t_0)\bar x+ t_0x_0\notin\qri (\Omega)$. Then Proposition~\ref{pro13} ensures the existence of $x^*\in X^*\setminus\{0\}$ with
\begin{equation}\label{tach1}
\left\langle x^*, x-\hat x\right\rangle \leq 0\;\text{ for all }\;x \in \Omega,
\end{equation}
\begin{equation}\label{tach2}
\left\langle x^*, \bar y- \hat x\right\rangle <0\;\text{ for some }\;\bar y\in \Omega.
\end{equation}
By the continuity of $x^*$ on $X$, we deduce from (\ref{tach1}) that
\begin{equation}\label{tach3}
\left\langle x^*, v \right\rangle \leq 0\;\text{ for all }\;v \in \overline{\cone}(\Omega-\hat x).
\end{equation}
On the other hand, the quasi-near convexity of $\Om$ yields the existence of a convex subset $C\subset X$ such that $\qri(C)\neq \emptyset$ and $C\subset \Omega \subset \overline{C}$. Thus $x_0\in \Bar{C}$ and by Proposition~\ref{pro2} we have $\ox\in \qri(\Bar{C})$. Using Proposition~\ref{pro3.12} gives us $\hat x\in \qri(\overline{C})$, i.e., $\overline{\cone}(\overline{C}-\hat x)$ is a subspace. Since $\overline{\cone}(\overline{C}-\hat x)=\overline{\cone}(\Omega-x_0)$, it follows from (\ref{tach3}) that
 $$\left\langle x^*, v \right\rangle = 0\text{ for all }v \in \overline{\cone}(\Omega-\hat x).$$
 This contradicts (\ref{tach2}) and completes the proof. $\h$

Now we give a characterization  of the {\em quasi-interior} of quasi-nearly convex sets.

\begin{proposition} \label{lem141}
Let $\Omega$ be a quasi-nearly convex set in $X$, and let $\bar x \in \Omega$. Then $\bar x \in \qi(\Omega)$ if and only if $N(\bar x;\Omega)=\{0\}$.
\end{proposition}
\textbf{Proof.}
Assume that $\bar x\in \qi (\Omega)$ and take any $x^* \in N(\bar x ; \Omega)$. Then the continuity of $x^*$ on $X$ ensures that $\langle x^*, v \rangle\leq 0$ for all $v \in \overline{\cone}(\Omega - \bar x)=X$, which yields $x^* = 0.$

Conversely, assume that $N(\bar x;\Omega)=\{0\}$. Arguing by contradiction, consider an arbitrary element $\bar v \in X$, and suppose that $\bar v \notin \overline{\cone}(\Omega-\bar x)$. Note that $\overline{\cone}(\Omega-\bar x)$ is convex as seen in the proof of Proposition \ref{pro2}(b). By \cite[Theorem~1.1.5]{zalinescu2002convex}, the sets $\{\bar v\}$ and $\overline{\cone}(\Omega-\bar x)$ can be  {\em strictly separated} by a closed hyperplane, i.e., there exists $x^* \in X^*\setminus\{0\}$ and $\alpha_1,\alpha_2 \in \mathbb{R}$ such that
$$
\left\langle x^*, v\right\rangle\leq\alpha_1<\alpha_2\leq\left\langle x^*, \bar v\right\rangle\;\text{ for all }\;v \in \overline{\cone}(\Omega-\bar x).
$$
Fix any $y \in \Omega$ and get $\left\langle x^*, y-\bar x\right\rangle\leq\alpha_1/\lambda$ as $\lambda>0$, which yields $\left\langle x^*, y-\bar x\right\rangle \leq 0$. Since $y$ was chosen arbitrarily in $\Omega$, we arrive at $x^* \in N(\bar x;\Omega)$, a contradiction completing the proof.
$\h$

The next proposition describes the separation of quasi-nearly convex sets via quasi-interiors.

\begin{proposition} \label{lem14}  Let $\Omega$ be a quasi-nearly convex set in $X$, and let $\bar x \in \Omega$. Then $\bar x \notin \qi(\Omega)$ if and only if the sets $\{\bar x\}$ and $\Omega$ can be separated by a closed hyperplane.
\end{proposition}
\textbf{Proof.}
By Proposition \ref{lem141}, we have $\bar{x} \notin \operatorname{qi}(\Omega)$ if and only if there is a nonzero element $x^*\in N(\bar{x}; \Omega)$. It follows  from the normal cone definition that $x^*\in N(\bar{x}; \Omega)$ if and only if
$$\left\langle x^*, y\right\rangle \leq\left\langle x^*, \bar{x}\right\rangle\;\text{ for all }\;y \in \Omega,$$
which completes the proof of the proposition. $\h$

By Definition~\ref{def1}(c,f) we have $\qi(\Omega) \subset\qri(\Omega)$. The next proposition shows that the equality holds therein if the quasi-interior of $\Om$ is nonempty.

\begin{proposition} \label{lem12}
Let $\Omega$ be a quasi-nearly convex set in $X$. If $\qi (\Omega)\neq \emptyset$, then $$\qi(\Omega)=\qri(\Omega).$$
\end{proposition}
\textbf{Proof.} Observe first that if $\qi(\Omega)\neq \emptyset$, then $N(0;\Omega -\Omega)=\{0\}$. Indeed, take any $x\in \qi(\Omega)$. Then $0\in \qi(\Omega-x)$, and since $\Omega -x\subset \Omega -\Omega$, we get $0\in \qi(\Omega-\Omega)$. Therefore, it follows from Proposition~\ref{lem141} that
\begin{equation}\label{key1.2}
N(0;\Omega -\Omega)=\{0\}.\end{equation}
To show that $\qri(\Omega)\subset\qi(\Omega)$, pick any $x\in \qri (\Omega)$ and $x^* \in N(0;\Omega -x)$. Then we get
\begin{equation}\label{key2}
\langle x^*, w-x\rangle\leq 0\;\text{ for all }\;w\in \Omega,
\end{equation}
which implies that $x^*\in N(x;\Omega)$. Since $x\in \qri (\Omega)$, it follows from Proposition~\ref{qrisubspace} that $N(x,\Omega)$ is a linear subspace. Thus we arrive at $-x^*\in N(x;\Omega)$, i.e.,
\begin{equation}\label{key2.1}
\langle x^*, x-w\rangle\leq 0\;\text{ for all }\;w\in \Omega.
\end{equation}
Summing up the two inequalities (\ref{key2}) and (\ref{key2.1}) gives us $$\langle x^*, w_1-w_2\rangle\leq 0\;\text{ whenever }\;w_1,w_2\in \Omega,$$
and so $x^*\in N(0; \Omega-\Omega)$. Then it follows from (\ref{key1.2}) that $x^*=0$. Since $x^* \in N(0;\Omega -x)$ was chosen arbitrarily, we get $N(0; \Omega-x) =\{0\}$. Applying finally Proposition~\ref{lem141} tells us that $x\in \qi(\Omega)$ and thus completes the proof of the proposition. $\h$

Now we ready to derive some {\em calculus rules} for quasi-relative interiors of quasi-nearly convex sets via continuous linear mappings.

\begin{theorem} \label{linearmap}
Let $T\colon X\to Y$ be a continuous linear mapping between LCTV spaces, and let $\Omega$ be a quasi-nearly convex set in $X$. Then the following assertions hold:
\begin{enumerate}
\item $T(\qri (\Omega))\subset \qri (T( (\Omega)).$
\item  If $T$ is injective and $T(\Omega)$ is quasi-regular, then $T(\qri(\Omega))= \qri (T(\Omega))$. The injectivity of $T$ is not required if $\Om$ is convex with $\qri(\Om)\neq\emptyset$.
\item $T(\Omega)$ and $T(\qri(\Omega))$ are quasi-nearly convex sets in $Y$.
\end{enumerate}
\end{theorem}
\textbf{Proof.} (a) Pick any $\bar x \in \qri(\Omega)$ and deduce from Proposition~\ref{qrisubspace} that $N(\bar x; \Omega)$ is a linear subspace. Take $y^* \in N(T(\bar x); T(\Omega))$ and get $\langle y^*, T(x) - T(\bar x)\rangle \leq 0$ for all $x \in \Omega$, which implies that $T^*y^* \in N(\bar x; \Omega)$. Since $N(\bar x;\Omega)$ is a subspace, we have that $-T^*y^* \in N(\bar x; \Omega)$. Equivalently, it holds that
$-y^* \in N(T(\bar x); T(\Omega))$, which implies that $N(T(\bar x); T(\Omega))$ is a linear subspace of $Y^*$. Hence Proposition~\ref{qrisubspace} tells us that $T(\bar x) \in \qri(T(\Omega))$.

(b) We have by the assumption that there is $\bar x \in \qri(\Omega)$, and thus $$\bar y: = T(\bar x) \in \qri(T(\Omega)).$$ Fix any $\hat y \in \qri(T(\Om))=\iri(T(\Omega))$ and deduce from Proposition~\ref{pro1} that there exists $y_0 \in T(\Omega)$ such that $\hat y=(1-t_0)\bar y+ t_0y_0\in T(\Omega)$ for some $t_0\in(0,1)$. Then we find $x_0,\hat x\in\Omega$ satisfying $y_0=T(x_0)$ and $\hat y = T(\hat x)$. Since $T$ is injective and $$T(\hat x)=\hat y=T((1-t_0)\bar x+t_0x_0),$$
we obtain $(1-t_0)\bar x+t_0x_0=\hat x\in\Omega$. Moreover, since $\bar x \in \qri(\Omega)$, Proposition~\ref{qriconvex} tells us that $\hat x \in \qri (\Omega)$, and hence $\hat y \in T(\qri(\Omega))$. The imposed assumptions and the proved result in (a) yield $$\qri (T(\Omega))\subset T(\qri(\Omega))\subset \qri (T(\Omega)).$$
(c) Since $\Om$ is quasi-nearly convex, we get from Definition~\ref{defqconvex} that there exists a convex set $C\subset X$ such that $\qri(C)\neq \emptyset$ and
$$C\subset \Om\subset \Bar{C}.$$
Then Proposition~\ref{pro2}(a) together with \cite[Proposition~2.12]{borwein1992partially} tells us that $\overline{\qri}(C)=\Bar{C}=\overline{\Omega}$. Thus
\begin{equation}\label{key13}
T(\qri(C))\subset T(\qri(\Omega))\subset T(\Omega)\subset T(\overline{\Omega})=T(\overline{\qri}(C))\subset \overline{T(\qri(C))},
\end{equation}
where the first inclusion is satisfied by Proposition~\ref{pro2}(b) and the last inclusion holds because $T$ is continuous. Note that $\qri(C)$ is convex by \cite[Lemma~2.9]{borwein1992partially}, and we deduce from the linearity of
$T$ that $T(\qri (C))$ is convex. Combining this with  (\ref{key13}) implies that
$T(\Omega)$ and $T(\qri(\Omega))$ are quasi-nearly convex.
$\h$

\section{Generalized Relative Interiors of Quasi-Nearly Convex Graphs}

The last section of the paper is devoted to deriving new results on quasi-relative interiors and quasi-interiors of graphs of set-valued mappings in the quasi-near convexity framework.

Given a {\em function} $f\colon X \rightarrow \overline{\mathbb{R}}$, we say that $f$ is {\em quasi-nearly convex} if $\epi(f)$ is a quasi-nearly convex set.  We also say that a {\em set-valued mapping} $F\colon X \rightrightarrows Y$ between LCTV spaces is {\em quasi-nearly convex} convex if its graph $\gph(F)$ is a quasi-nearly convex set in $X \times Y$.

The following proposition is useful to establish the major results of this section.

\begin{proposition}\label{prodomconvex}
Let $F\colon X \rightrightarrows Y$ is quasi-nearly convex set-valued mapping between LCTV spaces. Then both sets $\dom (F)$ and $\rge (F)$ are quasi-nearly convex sets. Consequently, if $f\colon X \rightarrow \overline{\mathbb{R}}$ is an extended-real-valued proper quasi-nearly convex function, then the sets $\dom (f)$ and $\rge (f)$ are quasi-nearly convex.
\end{proposition}
\textbf{Proof.}
We first see that $\dom (F)=\mathcal{P}(\gph (F))$, where $\mathcal{P}$ is the  continuous linear mapping given by
$$
\mathcal{P}(x, y):=x, \ \; (x,y) \in X \times Y.
$$
Then it follows from Theorem~\ref{linearmap}(c) that $\dom (F)$ is quasi-nearly convex. Similarly, $\rge (F)$ is quasi-convex because $\rge (F)=\mathcal{P}_1(\gph (F))$, where $\mathcal{P}_1$ is the continuous linear mapping defined by
$$
\mathcal{P}_1(x, y):=y, \ \; (x,y)\in X\times Y.
$$
Finally, suppose that $f$ is quasi-nearly convex. Then the epigraphical mapping $E_f$ defined in (\ref{epigraphical}) is quasi-nearly convex. Since $\dom (E_f)=\dom (f)$ and $\rge (E_f)=\rge (f)$, the sets $\dom (f)$ and $\rge (f)$ are both quasi-nearly convex, and we are done.
$\h$

Now we are ready for the first main result.

\begin{theorem} \label{qriqi}
Let $F\colon X \rightrightarrows Y$ be a quasi-nearly convex set-valued mapping. The following assertions are satisfied:
\begin{enumerate}
\item If $\graph(F)$ is quasi-regular, then we have
$$
\qri(\gph(F)) \subset\{(x, y) \in X \times Y \mid x \in \qri(\dom(F)), y \in \qri(F(x))\}.
$$
\item If $F(x)$ is quasi-nearly convex and $\qi(F(x))\neq \emptyset$ for every $x\in \dom (F)$, then $$\qri(\gph(F)) \supset \{(x, y) \in X \times Y \mid x \in \qri(\dom(F)), y \in \qri(F(x))\}.
$$
\end{enumerate}
Consequently, if $\graph(F)$ is quasi-regular and $\qi(F(x))\neq \emptyset$ for every $x\in \dom (F)$, then
$$
\qri(\gph(F)) = \{(x, y) \in X \times Y \mid x \in \qri(\dom(F)), y \in \qri(F(x))\}.
$$
\end{theorem}
\textbf{Proof.} (a) Observe first that $\dom (F)$ is quasi-nearly convex by Proposition \ref{prodomconvex}. We now fix any $(\bar{x}, \bar{y}) \in \qri ( \gph (F))$ and suppose on the contrary that $\bar{x} \notin\qri (\dom(F))$. Then it follows from Proposition~\ref{pro13} that the sets $\{\bar{x}\}$ and $\dom(F)$ can be properly separated by a closed hyperplane, i.e.,  there exist $x^* \in X^*\setminus\{0\}$ and $\hat x \in \dom(F)$ such that
$$
\left\langle x^*, x\right\rangle \leq\left\langle x^*, \bar{x}\right\rangle \text { for all } x \in \dom(F)
$$
and
$$\left\langle x^*, \hat x\right\rangle<\left\langle x^*, \bar{x}\right\rangle.$$
Thus, we have
$$
\left\langle\left(x^*, 0\right),(x, y)\right\rangle=\left\langle x^*, x\right\rangle \leq\left\langle x^*, \bar{x}\right\rangle=\left\langle\left(x^*, 0\right),(\bar{x}, \bar{y})\right\rangle
\text{ for all } (x, y) \in \gph(F)$$
and for each $\hat y \in F\left(\hat x\right)$,
$$
\left\langle\left(x^*, 0\right),\left(\hat x, \hat y\right)\right\rangle=\left\langle x^*, \hat x\right\rangle<\left\langle x^*, \bar{x}\right\rangle=\left\langle\left(x^*, 0\right),(\bar{x}, \bar{y})\right\rangle.
$$
This implies that the sets $\gph(F)$ and $\{(\bar{x}, \bar{y})\}$ can be properly separated. It follows from Proposition~\ref{pro13} that $(\bar{x}, \bar{y}) \notin \qri(\operatorname{gph}(F))$, a contradiction, which $\bar{x} \in \qri(\dom(F))$.
It remains to show that
$\bar{y} \in \qri(F(\bar{x}))$. Fix any $y \in F(\bar{x})$. Then by the quasi-regularity of $\gph(F)$ and Proposition \ref{pro1}, there exist $(x_0, y_0) \in \operatorname{gph}(F)$ and $t_0 \in(0,1)$ such that
$$
(\bar{x}, \bar{y})=(1-t_0)(\bar{x}, y)+t_0(x_0, y_0).
$$
This yields $x_0=\bar{x}$ and $\bar{y}=(1-t_0) y+t_0 y_0$ with $y_0 \in F(\bar{x})$.
 It follows therefore by
Proposition \ref{pro1} and Theorem \ref{2s}
that $\bar{y} \in \operatorname{iri}(F(\bar{x})) \subset \operatorname{qri}(F(\bar{x}))$, which completes the proof of (a).

(b) To verify the reverse inclusion in this assertion under the imposed assumptions that $F(x)$ is quasi-nearly convex and $\qi(F(x))\neq \emptyset$ for every $x\in\dom(F)$, we fix $\bar{x} \in \qri(\dom (F))$ and $\bar y \in \qri(F(\bar{x}))$. Arguing by contradiction, suppose that $(\bar{x}, \bar{y}) \notin \qri (\gph(F))$. Then it follows from Proposition~\ref{pro13} that there exist $\left(x^*, y^*\right) \in X^* \times Y^*\setminus\{(0,0)\}$ and $(\hat x,\hat y) \in X\times Y$ such that
\begin{equation} \label{46}
   \left\langle x^*, x\right\rangle+\left\langle y^*, y\right\rangle \leq\left\langle x^*, \bar{x}\right\rangle+\left\langle y^*, \bar{y}\right\rangle\;\text { whenever }\;x\in \dom(F\;)\text { and }\;y \in F(x)
\end{equation}
together with the strict inequality
\begin{equation} \label{f26}
\left\langle x^*, \hat x\right\rangle+\left\langle y^*, \hat y\right\rangle<\left\langle x^*, \bar{x}\right\rangle+\left\langle y^*, \bar{y}\right\rangle.
\end{equation}
We distinguish the two possible cases: (A) $y^*=0$ and (B) $y^*\neq 0$.

In case (A), we get from (\ref{46})  that
$$
\left\langle x^*, x\right\rangle \leq\left\langle x^*, \bar{x}\right\rangle \text { whenever } x\in \dom(F)
$$
and from (\ref{f26}) that
$$  \left\langle x^*, \hat x\right\rangle<\left\langle x^*, \bar{x}\right\rangle.$$
Then it follows from Proposition \ref{pro13} that $\bar{x} \notin \qri (\dom(F))$, a contradiction.

In case (B), letting $x=\bar{x}$ in (\ref{46}) gives us  $\left\langle y^*,  y\right\rangle \leq\left\langle y^*, \bar{y}\right\rangle$ for all $y \in F(\bar{x})$. Then Propositions~\ref{lem14} and \ref{lem12} tell us that $$\bar y\notin \qi(F(\bar x))=\qri(F(\bar x)).$$
This contradiction shows that $(\bar{x}, \bar{y}) \in \qri(\gph(F))$ and hence completes the proof of (B) and of the whole theorem.
$\h$

Another description of the quasi-relative interiors of graphs for quasi-nearly convex set-valued mapping formulated as follows.

\begin{theorem} \label{thmqrigraph}
Let $F\colon X \rightrightarrows Y$ be a quasi-nearly convex set-valued mapping. Then we have
$$\operatorname{qri}(\operatorname{gph}(F))\supset \{(x, y) \in X \times Y \mid x \in \operatorname{qri}(\operatorname{dom}(F)), y \in \operatorname{int}(F(x))\}.
$$
If in addition $\gph(F)$ is quasi-regular, $F(x)$ is convex and $\textrm{\rm int}(F(x))\neq \emptyset$ for every $x\in \dom (F)$, then
$$\operatorname{qri}(\operatorname{gph}(F)) = \{(x, y) \in X \times Y \mid x \in \operatorname{qri}(\operatorname{dom}(F)), y \in \operatorname{int}(F(x))\}.
$$
\end{theorem}
\textbf{Proof.}
 Pick any $(\bar{x}, \bar{y}) \in X \times Y$ with $\bar{x} \in \qri(\dom(F))$ and $\bar{y} \in \text{\rm int}(F(\bar{x}))$. We prove by contradiction, suppose that $(\bar{x}, \bar{y}) \notin \qri(\gph(F))$. By Proposition~\ref{pro13}, there exists $\left(x^*, y^*\right) \in X^* \times Y^*\setminus\{(0,0)\}$ such that
\begin{equation}\label{key223}
\langle x^*, x\rangle+\langle y^*, y\rangle \leq\langle x^*, \bar{x}\rangle+\langle y^*, \bar{y}\rangle \;\text { for all }\;x \in \dom(F)\;\text{ and }\;y \in F(x),
\end{equation}
and there exists $(\hat x,\hat y)\in\gph (F)$ for which
\begin{equation}\label{key233}
\langle x^*, \hat x\rangle+\langle y^*, \hat y\rangle <\langle x^*, \bar{x}\rangle+\langle y^*, \bar{y}\rangle.
\end{equation}
Letting $x=\bar{x}$ in (\ref{key223}), we obtain that
\begin{equation}\label{key244}
\langle y^*, y\rangle \leq\langle y^*, \bar{y}\rangle\;\text { whenever }\;y \in F(\bar{x}).
\end{equation}
Since $\bar{y} \in \sint(F(\bar{x}))$, there exists a symmetric neighborhood $V$ of the origin satisfying $V \subset \sint(F(\bar{x}))-\{\bar{y}\}$. It follows from (\ref{key244}) that $\left\langle y^*, v\right\rangle \leq 0$ and $\left\langle y^*,-v\right\rangle \leq 0$ for all $0 \neq v \in V$, and hence $y^*=0$ on $V$. Moreover, since $V$ is a symmetric neighborhood of $0\in X$, for every $0\neq y \in Y$ we find $0 \neq t \in \mathbb{R}$ such that $t y=v \in V$, and therefore $\left\langle y^*, y\right\rangle=\frac{1}{t}\left\langle y^*, v\right\rangle=0$. This implies that $y^*=0$ on $Y$, which implies by (\ref{key223}) and (\ref{key233}) that the sets $\{\bar{x}\}$ and $\dom(F)$ can be properly separated by a closed hyperplane. Then Proposition~\ref{pro13} tells us that $\bar{x} \notin \qri(\dom (F))$, a contradiction that justifies the inclusion
 $$\operatorname{qri}(\operatorname{gph}(F)) \supset \{(x, y) \in X \times Y \mid x \in \operatorname{qri}(\operatorname{dom}(F)), y \in \operatorname{int}(F(x))\}.
$$

To check the inclusion ``$\subset$" in the theorem under the additional assumptions made, we take any $(\bar{x}, \bar{y}) \in \qri(\gph(F))$ and suppose on the contrary that $\bar{x} \notin \qri(\dom(F))$. By Proposition~\ref{pro13}, there exist a nonzero function $x^*\in X^*$ and $\hat x\in \dom(F)$ such that
$$\langle x^*, x\rangle\leq \langle x^*,\bar x\rangle\;\text{ for all }\;x\in \dom (F)$$
and
$$\langle x^*,\hat x\rangle <\langle x^*,\bar x\rangle.$$
Then for all $(x,y)\in \gph (F)$, we get $$\left\langle (x^*,0), (x,y)\right\rangle= \langle x^*, x\rangle\leq \langle x^*,\bar x\rangle = \left\langle (x^*,0),(\bar x,\bar y)\right\rangle.$$
Taking any fixed $\hat y \in F(\hat x)$ gives us
$$\left\langle (x^*,0), (\hat x,\hat y)\right\rangle= \langle x^*, \hat x\rangle< \langle x^*,\bar x\rangle = \left\langle (x^*,0),(\bar x,\bar y)\right\rangle.$$
Thus, the sets $\{(\bar x,\bar y)\}$ and $\gph (F)$ can be properly separated. Applying Proposition~\ref{pro13} again, we arrive at the condition $(\bar{x}, \bar{y}) \notin \qri(\gph(F))$, which is a contradiction telling us that $\bar{x} \in \qri(\dom(F))$.

To proceed further, let us verify that $\bar{y} \in \text{int}(F(\bar{x}))$. Under the assumptions of the convexity of $F(\bar x)$ and the nonemptiness of $\text{\rm int}(F({\bar x}))$, we have by \cite[Theorem~2.12(b)]{borwein2003notions} that $\sint (F(\bar x))=\iri (F(\bar x))$. Therefore, it suffices to show $\bar{y} \in \text{iri}(F(\bar{x}))$. To justify this,
fix any $y \in F(\bar{x})$. Then the quasi-regularity of $\gph(F)$ and Proposition~\ref{pro1} give us $(x_0, y_0) \in \gph(F)$ such that
$$
(\bar{x}, \bar{y})=(1-t_0)(\bar{x}, y)+t_0(x_0, y_0)\;\text{ for some }\;t_0\in (0,1).
$$
This yields $x_0=\bar{x}$ and $\bar{y}=(1-t_0) y+t_0 y_0$ with $y_0 \in F(\bar{x})$. Applying Proposition~\ref{pro1} again tells that $\bar{y} \in \iri(F(\bar{x}))$, which completes the proof of the theorem.
$\h$

Next, we employ the above result to deduce representations of the quasi-relative interiors of {\em epigraphs} for extended-real-valued functions.

\begin{corollary}
 Let $f \colon X \to \oR$ be a proper quasi-nearly convex function. If $\epi (f)$ is quasi-regular, then we have
\begin{equation*}
\qri(\epi (f))=\big\{(x, \alpha)\in X\times \mathbb R\; \big |\; x\in \qri(\dom (f)), \; f(x)<\alpha\big\}.
\end{equation*}
\end{corollary}
\noindent {\bf Proof}. Consider the epigraphical mapping $E_f$ associated with $f$. Since $E_f(x)=[f(x),\infty)$ is convex and $\mbox{\rm int}(E_f(x))$ is nonempty every $x\in \dom (E_f)=\dom (f)$, the conclusion follows from Theorem~\ref{thmqrigraph}. $\h$

We end this section with the following inclusion of
quasi-interiors for graphs of quasi-nearly convex set-valued mappings.

\begin{theorem} \label{thmqi_graph}
Let $F\colon X \rightrightarrows Y$ be a quasi-nearly convex set-valued mapping between LCTV spaces. Then we have the inclusion
$$\operatorname{qi}(\operatorname{gph}(F)) \supset \{(x, y) \in X \times Y \mid x \in \operatorname{qi}(\operatorname{dom}(F)),\;y \in \operatorname{qi}(F(x))\}.
$$
\end{theorem}
\textbf{Proof.}  We proceed as in the proof of Theorem~\ref{qriqi}(b) with using now Proposition~\ref{lem14} instead of Proposition~\ref{pro13}.
 $\h$

\section*{Acknowledgements}
Research of Vo Si Trong Long is funded by University of Science,
VNU-HCM under grant number T2023-01.

\end{document}